\documentstyle[11pt,A4]{article}
\def\w{\omega}
\def\beq{\begin{equation}}
\def\eeq{\end{equation}}
\def \a{\alpha}
\def \b{\beta}
\def \g{\gamma}
\def \d{\delta}

\begin{document}

\begin{titlepage}

\title{A Link Invariant from Quantum Dilogarithm}

\author{R.M. Kashaev\thanks{On leave from
St. Petersburg Branch of the Steklov Mathematical Institute,
Fontanka 27, St. Petersburg 191011, RUSSIA}\\ \\
Laboratoire de Physique Th\'eorique
{\sc enslapp}\thanks{URA 14-36 du CNRS,
associ\'ee \`a l'E.N.S. de Lyon,
et \`a l'Universit\`e de Savoie}
\\
ENSLyon,
46 All\'ee d'Italie,\\
69007 Lyon, FRANCE\\
E-mail: {\sf rkashaev\@@enslapp.ens-lyon.fr}}

\date{April 1995}

\maketitle

\abstract{The link invariant, arising from the cyclic quantum dilogarithm
 via the particular $R$-matrix construction, is proved to coincide with the
invariant of triangulated links in $S^3$ introduced in \cite{K}. The obtained
 invariant, like Alexander-Conway polynomial, vanishes on disjoint union of
links. The $R$-matrix can be considered as the cyclic analog of the universal
$R$-matrix associated with $U_q(sl(2))$ algebra.}
\vskip 2cm

\rightline{{\small E}N{\large S}{\Large L}{\large A}P{\small P}-L-517/95}

\end{titlepage}

\section{Introduction}

Solutions to the Yang-Baxter equation (YBE) \cite{Y,B1} play important role
both in the theory of solvable two-dimensional models of statistical
mechanics and field theory (for a review see \cite{B2,F}), and in the theory
of braid groups and links \cite{Jo1}. Algebraically YBE is deeply related with
quantum groups \cite{D1}.

Recently, Rogers' dilogarithm identity \cite{R} attracted much attention in
connection with solvable models and 3d mathematics, see \cite{DS} and
references
therein. In \cite{FK} quantum generalization of the dilogarithm function has
been suggested. Remarkably, quantized dilogarithm identity, like YBE, is also
connected with quantum groups \cite{K,K1}.

In \cite{K} the invariant of triangulated links in three-manifolds has been
defined. It was conjectured that the choice of the triangulation is
inessential.
In this paper we prove this conjecture for the case of three-sphere. Namely,
using particular $R$-matrix, we define in the standard way the link invariant,
 and then show that it coincides with the three-dimensional construction of
\cite{K}.

In Sect.\ref{sec1} the $R$-matrix is defined. The corresponding link invariant
is constructed in Sect.\ref{sec2}. In Sect.\ref{sec3} the three-dimensional
description is developed and the equivalence with construction of \cite{K}
is established.

\section{The $R$-Matrix}
\label{sec1}

First introduce the necessary notations. Let $\w$ be some primitive $N$-th root
of unity for some $N\ge 2$.
For any non-negative integer $n<N$ define the symbol:
\beq
(\w)_n=\cases{1, &if $n=0$;\cr
              \prod_{j=1}^n(1-\w^j), &$0<n<N$,\cr}
\eeq
Under complex conjugation we have the following properties:
\beq
(\w)_n^*=(\w^*)_n,\quad \w^*\w=1.
\eeq
It will be convenient to use the ``theta''-function, defined on integers:
\beq
\theta(n)=\cases{1, &$0\le n< N$;\cr
                 0,&otherwise,\cr}\quad n\in Z.
\eeq
For any $n\in Z$ denote by $[n]$ the coresponding
integer residue modulo $N$, lying in the interval $\{0,\ldots,N-1\}$:
\beq
[n]=n\pmod N,\quad 0\le [n]<N,\quad n,[n]\in Z.\label{brac}
\eeq
The evident property of this function is its' periodicity:
\beq
[n+N]=[n],\quad n\in Z. \label{periodicity}
\eeq
It should be stressed however that the ``bracketing'' is a map from $Z$ to
 $Z$ rather than to $Z_N$, the ring of residues modulo $N$. Besides,
 it does not preserve the addition operation:
\beq
[m]+[n]\ge N\Rightarrow [m]+[n]\ne [m+n].
\eeq
Define a function with four integer arguments:
\beq
W(k,l,m,n)=N\w^{(l+m)(m+n)}{\theta([k]+[m])\theta([l]+[n])\over
(\w)_{[k]}(\w^*)_{[l]}(\w)_{[m]}(\w^*)_{[n]}},\quad k,l,m,n\in Z.\label{4W}
\eeq
Using these definitions consider the following indexed set of numbers:
\beq
R(i,j,k,l|a,b,c,d)=\w^{a-k-j}W(j-i-a,i-l-d,l-k-c,k-j-b),
\label{R}
\eeq
where all indices are integers, and $a,b,c,d$ being constrained to have
the unit sum modulo $N$:
\beq
a+b+c+d=1\pmod N.
\eeq
There is a symmetry group of these symbols generated by relations:
\beq
R(i,j,k,l|a,b,c,d)=R(k,l,i,j|c,d,a,b),             \label{symmetry1}
\eeq
and
\beq
R(i,j,k,l|a,b,c,d)=R(-j,-i,-l,-k|a,d,c,b)\w^{-i-j-k-l},\label{symmetry2}
\eeq
which have a natural geometrical interpretation. Indeed, consider the
diagram given by a projection of
two non-coplanar segments in $R^3$ to some plane, the images of the segments
having an intersection point. Place indices $j$ and $l$ on two ends of the
overcrossing segment, and indices $i$ and $k$, on the ends of the other
segment, choosing among two possibilities that, where indices $i,j,k,l$
go round the intersection point in the counter-clockwise direction. Place
indices $a,b,c,d$ on four sectors also in the counter-clockwise order,
starting from the sector bounded by segments with the end labels $i$ and $j$.
Associate with this configuration the ``Boltzmann'' weight (\ref{R}). The above
mentioned symmetry group acquires now an interpretation of the symmetry
group of the diagram. Namely, relation (\ref{symmetry1}) corresponds to the
rotation by $\pi$ aroung the intersection point, while (\ref{symmetry2})
corresponds to looking at the projection image from the opposite side of the
plane.

The next important property of the $R$-symbols is the Yang-Baxter identity,
which they satisfy. Using periodicity property (\ref{periodicity}), confine
the range, where the indices take their values, to the finite set
$\{0,\ldots,N-1\}$. Consider an $N^2$-by-$N^2$ $R$-matrix,
described by the matrix elements:
\beq
\langle i,j|R|k,l\rangle=R(i,j,k,l|1,0,0,0)\w^{k+l}.\label{R-matrix}
\eeq
Then, the following YBE holds:
\beq
R_{12}R_{13}R_{23}=R_{23}R_{13}R_{12},\label{YBE}
\eeq
where the standard matrix notations are used:
\beq
R_{12}|i,j,k\rangle=\sum_{i',j'=0}^{N-1}|i',j',k\rangle
\langle i',j'|R|i,j\rangle,
\eeq
and similarly for the other matrices. Later we will argue that the $R$-matrix
(\ref{R-matrix}) is in fact a cyclic analog of the universal $R$-matrix
associated with the quantum group $U_q(sl(2))$ (see the remark in the end of
Sect.\ref{sec3}).

The inverse $R$-matrix is given in terms of the same $R$-symbols:
\beq
\langle i,l|\overline R|k,j\rangle=R(i,j,k,l|,0,0,0,1)\w^{j+k}.
\eeq
So, in matrix notations we have
\beq
R_{12}\overline R_{12}=1.
\eeq
Thus, taking into account above diagrammatic interpretation as well as the
symmetry property (\ref{symmetry1}), we have already
realizations of two Reidemeister moves, corresponding to the regular isotopy
of links. The following relation provides us with the third Reidemeister move,
necessary for the ambient isotopy:
\beq
\sum_{j=0}^{N-1}\w^jR(i,j,j,l|0,1,0,0)=\w^{-i}\delta_{i,l}.
\eeq
In the next section, using these results, we define the ``state sum'',
leading to an ambient isotopy link invariant.

\section{Invariants of Links in Three-Sphere}
\label{sec2}

\setcounter{equation}{0}

Let link $L$ in $S^3$ is given by a two-strand tangle. Consider some
non-singular (i.e. with only simple intersection points) planar projection
of this tangle. Denote by $V,E,F$ the sets of vertices, edges, and faces,
 respectively, in the planar graph. We will write $v\in f$, for some vertex
$v$ and face $f$ if $v$ lies in the boundary of $f$. There are two
distinguished elements $e_1$
 and $e_2$ in $E$, corresponding to the two strands of the tangle,
and one distinguished element $f_0$ in $F$, corresponding to the outer
 region in the plane. Define the following two maps:
\beq
 j\colon E\to Z_N,\quad a\colon V\times F\to Z_N,\label{ja-maps}
\eeq
where we identified $Z_N$ with the set $\{0,\ldots,N-1\}$. Impose the
restrictions on $j$:
\beq
j(e_1)=j(e_2)=0,\label{boundary}
\eeq
and on $a$:
\beq
a(v,f_0)=0;\quad v\not\in f\Rightarrow a(v,f)=0,\quad  v\in V; f\in F,
\eeq
\beq
\sum_{f\in F}a(v,f)=\sum_{v\in V}a(v,f)=1,\quad f\ne f_0.
\eeq
Using these maps, associate with each vertex the Boltzmann weight
according to Sect.~\ref{sec1}. This gives us one more map:
\beq
r_{j,a}\colon V\to C.
\eeq
Consider now the ``partition'' function:
\beq
\langle L\rangle=\sum_j \prod_{v\in V}r_{j,a}(v)\prod_{e\in E}\w^{j(e)}.
\label{inv}
\eeq

\noindent
\newtheorem{theorem}{Theorem}
\begin{theorem}
(1) The quantity $\langle L\rangle^N$ is an ambient isotopy invariant of
the link $L$ in three-sphere; (2) for odd $N$ the following equality holds:
\beq
\langle L\rangle^N=Q(S^3,L),
\eeq
where $Q(M,L)$ is the invariant of triangulated link $L$ in three-manifold $M$
defined in~\cite{K}.\label{theo}
\end{theorem}
The first part of the theorem is a consequence of the identities, satisfied by
the $R$-symbols (see Sect.~\ref{sec1}). There is a freedom in choosing the
 second map in (\ref{ja-maps}). In fact, different choices lead to one and
 the same partition function (\ref{inv}) up to $N$-th roots of unity. Note
that the face indices are nothing else but the angle dependent data,
corresponding to ``enhanced'' vertex models \cite{BKW,Jo,W}. This becomes
evident after Fourier transformation over the edge indices. The face
index - angle correspondence is described by the formula:
\beq
\phi=\pi-2\pi a,
\eeq
where $\phi$ is the angle between the adjoining straight edges in the
 diagram, corresponding to $R$-symbol (see Sect.~\ref{sec1}), and $a$
 is the face index associated with the sector bounded by these edges.

The proof of the second part of the theorem is given in the next section.

Note that the two strand tangle representation of links has been used in
\cite{Kauf}
to construct the state sum for Alexander-Conway polynomial \cite{Alex,Conw},
because the latter is zero for disjoint union of links. This is also true for
 the invariant (\ref{inv}):
\beq
\langle L_1\sqcup L_2\rangle=0.\label{vanish}
\eeq
In the language of quantum group theory, property (\ref{vanish}) reflects
the vanishing quantum dimension. In the three-dimensional definition given
in \cite{K} ( see also the next section) no tangle description is needed.

\section{Three-Dimensional Description}
\label{sec3}
\setcounter{equation}{0}

In this section we establish the connection of the link invariant constructed
in Sect.~\ref{sec2} with the invariant of triangulated links of \cite{K},
related in turn to cyclic quantum dilogarithm \cite{FK}.
First, remind the necessary formulas from \cite{K} in slightly different
notations.

We will assume that the integer $N$ is odd, and the square root of $\w$ will
be thought to be taken as an $N$-th root of unity:
\beq
\w^{1/2}=\w^{(N+1)/2}.
\eeq
For complex $x,y,z$, satisfying the Fermat equation
\beq
x^N+y^N=z^N,\label{Fermat}
\eeq
and integers $m,n$ define the function
\beq
w(x,y,z|m,n)=w(x,y,z|m-n)\w^{n^2/2},
\eeq
where
\beq
w(x,y,z|n)=\prod_{j=1}^n{y\over z-x\w^j}.
\eeq
The Fermat equation (\ref{Fermat}) ensures the periodicity relation:
\beq
w(x,y,z|n+N)=w(x,y,z|n).
\eeq
This latter property enables us to consider the finite set $Z_N$ instead of
$Z$.

For any non-coinciding complex numbers $z_0,z_1,z_2,z_3$, and
$\a,\b,\g,\d,a,c\in Z_N$ define two functions:
\beq
T^{0123}(a,c)_{\a,\b}^{\g,\d}=
\rho_{0123}\w^{c(\g-\a)+\a\d+ac/2}w(x_{03}x_{12},x_{01}x_{23},
x_{02}x_{13}|\g-a,\a)\delta_{\b,\g+\d},\label{T}
\eeq
and
\beq
\overline T^{0123}(a,c)^{\a,\b}_{\g,\d}=
{\overline\rho_{0123}\w^{c(\g-\a)-\a\d-ac/2}\delta_{\b,\g+\d}\over
w(x_{03}x_{12}/\w,x_{01}x_{23},
x_{02}x_{13}|\g+a,\a)},\label{Tbar}
\eeq
where
\beq
x_{ij}=(z_i-z_j)^{1/N},\quad i,j=0,1,2,3,
\eeq
with $N$-th root being chosen to be real for real difference $z_i-z_j$;
scalar functions $\rho_{0123}$ and $\overline\rho_{0123}$ depend on only $z$'s
and their explicit form can be found in \cite{K}.
 These symbols
can be associated with tetrahedrons in $R^3$ as is described below.

In what follows for a polyhedron $X$, considered as a collection
of vertices, edges, triangles, and tetrahedrons,
 we will denote by $\Lambda_i(X)$, $i=0,1,2,3$ the sets of
simplices of corresponding dimension.
Consider a topological \footnote{
the term ``topological'' here means that edges and faces of the
tetrahedron can be curved}
 tetrahedron $t$ in $R^3$.
Order the vertices by fixing a bijective map
\beq
u\colon \{0,1,2,3\}\ni i\mapsto u_i\in\Lambda_0(t).
\eeq
 Put an arrow on each edge, pointing from a ``larger'' vertex
(with respect to above ordering) to a ``smaller'' one.
Each face also gets its' own orientation.
Let $u_l$ be the top of the tetrahedron.
Let us look from it down at the vertices $u_i,u_j,u_k$, where
\beq
\{i,j,k\}=\{0,1,2,3\}\setminus \{l\},\quad i<j<k.
\eeq
We will see two possible views: either $u_i,u_j,u_k$, in the order which they
are written, go round in the counter-clockwise direction
(the ``right'' orientation) or, in the
clockwise one (the ``left'' orientation).
The tetrahedron
 itself has two possible orientations in the following sense. Call the
tetrahedron ``right (left)''-oriented if the face $u_0u_1u_2$ has the
right (left) orientation in the above sense.

Introduce three maps:
\beq
s\colon \Lambda_0(t)\to C,\quad c\colon \Lambda_1(t)\to Z_N,
\quad\alpha\colon \Lambda_2(t)\to Z_N,\label{charge}
\eeq
where $s$ is injective, and $c$ satisfies the following relations:
\beq
\sum_{e\ni v}c(e)=1/2,\quad v\in\Lambda_0(t).
\eeq
Let $c_{ij}=c(u_iu_j)$ ($u_iu_j$ is the edge having ends $u_i$ and $u_j$),
$\alpha_i=\alpha(u_ju_ku_l)$, $\{j,k,l\}=\{0,1,2,3\}\setminus\{i\}$, and
$z_i=s(u_i)$.
Define the symbol associated with the tetrahedron $t$:
\beq
t_u(s,c,\alpha)=\cases{ T^{0123}(c_{01},c_{12})
_{\alpha_3,\alpha_1}^{\alpha_2,\alpha_0},& right orientation;\cr
\overline T^{0123}(c_{01},c_{12})
^{\alpha_3,\alpha_1}_{\alpha_2,\alpha_0},& left orientation.\cr}\label{tetr}
\eeq
The variables associated with the second map in (\ref{charge}) will be called
$Z_N$-charges. Using symbols (\ref{tetr}), one can associate a ``partition''
function to any finitely triangulated three-manifold with a triangulated tangle
passing through the all interior 0-simplices. The construction goes as follows.

Let $M$ be a finite triangulation of an oriented 3-dimensional
manifold with boundary. Fix a subset of 1-simplices
$T\subset\Lambda_1(M\setminus\partial M)$ in such
a way that any interior 0-simplex belongs to exactly two elements from $T$,
so the latter
is some triangulated tangle in $M$, passing through the all interior vertices.
 Denote
$I=\{0,1,\ldots,K-1\}$, where $K$ is the number of vertices in $M$, and fix
the following maps:
\beq
u\colon I\to\Lambda_0(M),\quad s\colon\Lambda_0(M)\to C,\quad
c_T\colon\Lambda_3(M)\times\Lambda_1(M)\to Z_N,\quad
\alpha\colon\Lambda_2(M)\to Z_N,\label{allmaps}
\eeq
where $u$ is bijective, $s$, injective, and $c_T$ satisfies the restrictions:
\beq
e\not\in\Lambda_1(t)\Rightarrow c_T(t,e)=0;
\quad\sum_{e\ni v}c_T(t,e)=1/2,\quad v\in t\in\Lambda_3(M);
\eeq
\beq
\sum_{t\in\Lambda_3(M)}c_T(t,e)=\cases{0,&$e\in T$;\cr
				       1,&$e\in
\Lambda_1(M\setminus\partial M)\setminus T$.\cr}
\eeq
Denote the total $Z_N$-charges on $\Lambda_1(\partial M)$ by
\beq
\partial c_{tot}\colon \Lambda_1(\partial M)\ni e\mapsto
\sum_{t\in\Lambda_3(M)}c_T(t,e)\in Z_N,
\eeq
and also the restrictions of the remained maps in (\ref{allmaps}) by
\beq
\partial u\colon\Lambda_0(\partial M)\to I,\quad
\partial s\colon\Lambda_0(\partial M)\to C,\quad
\partial\alpha\colon\Lambda_2(\partial M)\to Z_N.
\eeq
Then, the partition function reads:
\beq
M_{\partial u}(T,\partial s,\partial c_{tot},\partial\alpha)=N^{-K_0}
\sum_{\alpha,\partial\a=\mathrm{fixed}}
\prod_{t\in\Lambda_3(M)}
t_u(s,c_T,\alpha)\prod_{e\in\Lambda_1(M\setminus\partial M)\setminus T}
\langle s(\partial e)
\rangle^{-1},\label{partfun}
\eeq
where
\beq
\langle s(\partial e)\rangle=(s(u_j)-s(u_i))^{(N-1)/N},\quad e=u_iu_j;
\eeq
and $K_0$ being the number of interior 0-simplices. As is reflected in the
left hand side of (\ref{partfun}), the partition function depends (up to $N$-th
roots of unity) only
on the boundary data as well as on the topological class of the triangulated
tangle (see \cite{K} for the definition). This is the consequence of the
relations the symbols (\ref{tetr}) satisfy, which are described in \cite{K}.
When $\partial M=\emptyset$ partition function (\ref{partfun}) multiplied by
$N^2$ reduces to the invariant of triangulated links \cite{K}.

To make contact with the results of Sect.~\ref{sec2}, we have just to specify
the triangulated 3-manifold with tangle, corresponding to $R$-symbol (\ref{R}).
The most economic way is to consider singular triangulated manifolds \cite{TV}.

Let $u_1,u_2,u_3,u_4$ be four points in $R^3$, specifying the
right-oriented tetrahedron $u_1u_2u_3u_4$. Glue to four faces of the latter
 four
other tetrahedrons $$u_1u_2u_4u_5,\quad u_2u_3u_4u_{5'},\quad u_0u_1u_2u_3,
\quad u_{0'}u_1u_3u_4.$$ We have added four new vertices
$u_0,u_{0'},u_5,u_{5'}$. Among the latter identify the ``primed'' vertices with
the corresponding ``non-primed'' ones:
\beq
u_0=u_{0'},\quad u_5=u_{5'}.
\eeq
While doing this, we glue also the added tetrahedrons pairwise along the
corresponding faces. As a result we get a triangulated octahedron
$u_0u_1u_2u_3u_4u_5$. The next step consists in attaching four more
tetrahedrons, $$u_0u_1u_2u_5,\quad u_0u_2u_3u_5,\quad u_0u_3u_4u_5,
\quad u_0u_1u_4u_5,$$
to the octahedron, each one being glued along two faces. For example,
tetrahedron $u_0u_1u_2u_5$ is attached along faces $u_0u_1u_2$ and $u_1u_2u_5$,
and similarly the others. Thus we obtain singular triangulated three-ball $B$,
where among others four different 1-simplices have coinciding ends $u_0$
 and $u_5$. We will distinguish them by specifying the tetrahedron to
 which they belong. For example,
\beq
[u_0u_5]_{12}\in \Lambda_1(u_0u_1u_2u_5),
\eeq
and similarly for the others. Besides, $\Lambda_2(\partial B)$ consists of four
pairs of 2-simplices, each being given by one and the same triple of vertices:
\beq
\Lambda_2(\partial B)=\{[u_0u_iu_5]^{\pm},\quad i=1,2,3,4\},
\eeq
where we distinguish elements within each pair by their orientation ($+$ ($-$)
 corresponds
to the right (left) orientation).
To write down the partition function, we have to specify the tangle. Let us
choose it as follows:
\beq
T=\{u_1u_3,u_2u_4\}.
\eeq
Now, providing $B$ with additional data (\ref{allmaps}), one can calculate
 the partition function (\ref{partfun}). The result reads:
\begin{eqnarray}
\lefteqn{B_{\partial u}(T,\partial s,\partial c_{tot},
\partial\alpha)}\qquad\nonumber\\
&&=(s_{05})^{N-1}
\sum_{i_1,i_2,i_3,i_4=0}^{N-1}R(i_1,i_2,i_3,i_4|c_{12},c_{23},c_{34},c_{14})
\prod_{m=1}^4\psi_{m,i_m}^{\a_m^--\a_m^++c_m},\label{octahed}
\end{eqnarray}
where
\beq
 \a_m^\pm=\alpha([u_0u_mu_5]^{\pm}),\quad s_{ij}=(s(u_i)-s(u_j))^{1/N}
\eeq
\beq
c_m=\partial c_{tot}(u_mu_5),\quad c_{ij}=\partial c_{tot}([u_0u_5]_{ij}),
\eeq
and
\beq
\psi_{m,j}^k=\w^{jk}s_{0m}^{[k-1]}s_{m5}^{[-k]},
\eeq
see also (\ref{brac}). The role of $\psi$-functions is to intertwine the
 indices of different nature. Note that the $Z_N$-charges are related with the
 angle-dependent data for the enhanced $R$-matrix.

Formula (\ref{octahed}) enables us to prove the second part of Theorem
\ref{theo}. Indeed, let
we are given a link, represented by a two-strand tangle. Then, using
correspondence (\ref{octahed}), we associate to the latter partition function
of a singular triangulated three-ball with four boundary 0-simplices, the
 tangle
by its strands being attached to two of them. To get the invariant in $S^3$,
we have to glue the ball with another (empty) one, making the tangle to pass
through the other two vertices and converting it into the link under
 consideration. This last process is equivalent to choosing boundary conditions
(\ref{boundary}).

One can prove also the YBE (\ref{YBE}) via consideration the triangulated
 three-balls, corresponding to two sides of YBE, and transforming them to
each other by a sequence of elementary moves defined in \cite{K}.

A remark is in order. In \cite{K1} it has been shown that the universal
$R$-matrix in Drinfeld double can be represented as a product of solutions
to the constant pentagon relation, the canonical element in Heisenberg double
being such solution. This construction can be extended also to the case of
 non-constant solutions to the pentagon relation. In particular, using cyclic
 quantum dilogarithm of \cite{FK} and specializing the parameters, one can
derive $R$-matrix
(\ref{R-matrix}). The derivation given in this section is essentially
 equivalent
to the latter. On the other hand, the cyclic quantum dilogarithm can be
 obtained from the quantum dilogarithm with $q$ non-root of unity
\footnote{in this case the quantum dilogarithm is precisely the canonical
 element in Heisenberg double of $BU_q(sl(2))$,  the Borel sub-algebra of
 $U_q(sl(2))$ \cite{K}} in the limit, where $q$ approaches a root of unity
\cite{BR}. Thus, it is natural to consider $R$-matrix (\ref{R-matrix}) as a
 cyclic counterpart of the universal $R$-matrix in $U_q(sl(2))$.

\section{Summary}

According to \cite{K}, cyclic quantum dilogarithm leads to invariants of
 triangulated links in triangulated three-manifolds. It was conjectured that
these invariants are independent of triangulations used. The $R$-matrix
derivation, presented in this paper (Sect.\ref{sec2}), proves the conjecture
for the case of three-sphere. The essential part of the proof is the
realization
of the $R$-matrix (\ref{R-matrix}) as a three-dimensional partition function
(\ref{partfun}), (\ref{octahed}) for the triangulated three-ball with
properly chosen tangle in it.

The $R$-matrix (\ref{R-matrix}) can
be considered as a cyclic counterpart of the universal $R$-matrix for the
algebra $U_q(sl(2))$. The factorized structure of the matrix elements
(\ref{R}), (\ref{4W}) probably reflects the fact that this $R$-matrix
can be obtained in some limit from the chiral Potts model $R$-matrix
\cite{BPA}.

{\bf Acknowledgements.} The author thanks L.D. Faddeev for his encouragement
in this work, and  P. Degiovanni and L. Freidel for discussions.

\end{document}